\documentclass{amsart}
\newtheorem{thm}{Theorem}[section]

\newtheorem{prop}[thm]{Proposition}
\newtheorem{cor}[thm]{Corollary}

\theoremstyle{remark}
\newtheorem{rem}[thm]{Remark}

\numberwithin{equation}{section}
\newcommand{\cL}{\mathcal{L}}
\newcommand{\N}{\mathbb{N}}
\newcommand{\BA}{\ensuremath{\mathbb{A}}}
\newcommand{\GG}{\ensuremath{\mathbb{G}}}
\newcommand{\HH}{\ensuremath{\mathbb{H}}}
\DeclareMathOperator{\az}{alph}
\DeclareMathOperator{\ord}{ord}
\begin{document}
\title[Elementary equivalence of graph products]{Elementary equivalence of right-angled Coxeter groups and graph products of finite abelian groups}

\author[M. Casals-Ruiz]{Montserrat Casals-Ruiz}\thanks{The first author is supported by Programa de Formaci\'{o}n
de Investigadores del Departamento de Educaci\'{o}n, Universidades e Investigaci\'{o}n del Gobierno Vasco}
\address{Montserrat Casals-Ruiz,
Department of Mathematics and Statistics, McGill University,
805 Sherbrooke St. West, Montreal,
Quebec H3A 2K6, Canada}
\email{casalsruiz@math.mcgill.ca}

\author[I. Kazachkov]{Ilya Kazachkov}
\address{Ilya V. Kazachkov,
Department of Mathematics and Statistics, McGill University,
805 Sherbrooke St. West, Montreal,
Quebec H3A 2K6, Canada}
\email{kazachkov@math.mcgill.ca}
\thanks{The second author is supported by the Bourse d'excellence de l'Institut des sciences Math\'{e}matiques}

\author[V. Remeslennikov]{Vladimir Remeslennikov}
\address{Vladimir Remeslennikov,
Omsk Branch of Institute of Mathematics (SB RAS),
13 Pevtsova St., Omsk, 644099, Russia}
\email{remesl@iitam.omsk.net.ru}
\thanks{The third author is supported by the RFBR grant \#08-01-00067}

%20E06 Free products, free products with amalgamation, Higman-Neumann-Neumann extensions, and generalizations
%20F55 Reflection and Coxeter groups
%03C07 Basic properties of first-order languages and structures
%20E99 structure and classification of infinite groups: None of the above, but in this section

\begin{abstract}
We show that graph products of finite abelian groups are elementarily equivalent if and only if they are $\exists\forall$-equivalent if and only if they are isomorphic. In particular, two right-angled Coxeter groups are elementarily equivalent if and only if they are isomorphic.
\end{abstract}
\maketitle

The notion of elementary equivalence is fundamental in model theory. One of the most natural problems about elementary equivalence, is given a class of algebraic systems, to understand which systems in this class are elementarily equivalent, i.e. to classify algebraic systems up to elementary properties. This problem in the case of groups is usually rather hard and there are only few examples known when this problem has a satisfactory solution.

For the class of abelian groups, this is a well-known result of W.~Szmielew (1955). For the class of ordered abelian groups this problem was studied by A.~Robinson and E.~Zakon (1960), M.~Kargapolov (1963) and Yu.~Gurevich (1964).

A.~Malcev (1961) solved this problem for classical linear groups. His approach was generalised by E.~Bunina and A.~Mikahl\"{e}v to other linear, algebraic and Chevalley groups. The problem of classifying linear groups over integers up to elementary properties was studied by V.~Durnev (1995).

For certain nilpotent groups this problem was studied by O.~Belegradek, R.~Deborah, A.~Miasnikov, F.~Oger, V.~Remelsennikov. For certain free operator groups, this problem was solved by A.~Miasnikov and V.~Remeslennikov (1987).

Around 1945 it was conjectured by A.~Tarski, that elementary theories of free non-abelian groups of different rank coincide. This conjecture is now known as Tarski's problem and has recently been solved by O.~Kharlampovich and A.~Miasnikov (2006), and, independently, by Z.~Sela (2006). The classification of torsion-free hyperbolic groups up to elementary properties has been recently announced by Z.~Sela.

In this paper we address the classification of graph products of finite abelian groups up to elementary properties. We prove that two graph products of finite abelian groups are elementarily equivalent if and only if they are $\exists\forall$-equivalent if and only if they are isomorphic. In particular, two right-angled Coxeter groups are elementarily equivalent if and only if they are isomorphic.

\section{Preliminaries}
\subsection{Graph products of groups}
The idea of a graph product, introduced in \cite{Green}, is a generalisation of the concept of a partially commutative group. Let $G_1,\dots, G_k$, $G_i=\langle X_i\mid R_i\rangle$, $i=1,\dots, k$ be groups.  Let $\mathcal{G}_0=(V(\mathcal{G}_0), E(\mathcal{G}_0))$ be a finite, undirected, simplicial graph, $V(\mathcal{G}_0)=\{v_1,\dots,v_k\}$.

We now define the notion of a graph product of groups and state the necessary preliminary results. We refer the reader to \cite{Green} and \cite{Goda} for proofs.

A graph product $\GG=\GG_{\mathcal{G}}(A_1,\dots,A_k)$ of the groups $A_1,\dots, A_k$ with respect to the graph $\mathcal{G}_0$, is a group with a presentation of the form
$$
\langle X_1,\dots, X_k\mid R_1,\dots, R_k, \mathcal{R}\rangle,
$$
where $\mathcal{R}=\{[X_i,X_j]\mid \hbox{$v_i$ and $v_j$ are adjacent in $\mathcal{G}_0$}\}$.

In this paper we work with a graph product $\GG$ of finite abelian groups $A_1,\dots, A_k$. Denote the order of the group $A_i$ by $\ord(A_i)$. An important case of the graph product of finite abelian groups is the case when $\ord(A_i)=2$ for all $i$, in which case the group $\GG$ is a right-angled Coxeter group. These groups are widely studied, see \cite{Davis} and references there.

It is convenient to encode the graph product of finite abelian groups $\GG$ by a \emph{marked graph} $\mathcal{G}$. The marked graph $\mathcal{G}$ is obtained from the graph $\mathcal{G}_0$ as follows. Every vertex group $A_i$ of $\GG$ is the direct product of finite cyclic groups, whose orders are powers of primes, $$
A_i=G_{i,1}\times\dots\times G_{i,r_i},
$$
where $\ord(G_{i,j})=p_{i,j}^{n_{i,j}}$, $n_{i,j}\in\N$ and $p_{i,j}$ is a prime. The above decomposition of $A_i$ is unique up to a permutation of the factors. We replace the vertex $v_i$ in $\mathcal{G}_0$ by $r_i$ vertices $v_{i,1},\dots,v_{i,r_i}$. There is an edge in $\mathcal{G}$ between $v_{i,j_1}$ and $v_{i,j_2}$ for every $i$, $j_1$ and $j_2$, $1\le i\le k$, $1\le j_1<j_2\le r_i$, i.e. the full subgraph of $\mathcal{G}$ on $v_{i,1},\dots,v_{i,r_i}$ is a complete graph. Furthermore, if there is an edge  between $v_i$ and $v_j$ in $\mathcal{G}_0$, then there is an edge between $v_{i,l_1}$ and $v_{j,l_2}$ in $\mathcal{G}$, for all $1\le l_1\le r_i$, $1\le l_2\le r_j$.

Note that the graph product of groups $G_{1,1},\dots, G_{1,k_1},G_{2,1},\dots, G_{k,r_k}$ with the underlying graph $\mathcal{G}$ is isomorphic to $\GG$. The vertex groups $G_{i,j}$ of this graph product are finite, cyclic and directly indecomposable. We mark the vertices $v_{i,j}$ of the graph $\mathcal{G}$ by the orders $p_{i,j}^{n_{i,j}}$ of the corresponding vertex group $G_{i,j}$, where the $p_{i,j}$'s are primes and $n_{i,j}\in \N$.

We further assume that every graph product of abelian groups $\GG$ is given by the marked graph $\mathcal{G}$ constructed above. In other words, we treat $\GG$ as a graph product of finite, cyclic, directly indecomposable groups. Elisabeth Green proved in \cite{Green} that if a group can be represented as a graph product of cyclic groups of prime order, then this representation is unique. This result was extended to primary cyclic groups by D.~Radcliffe in his PhD Thesis. In \cite{Radcliffe}, Radcliffe showed that if a group can be represented as a graph product of directly indecomposable finite groups, then this representation is unique.

Let  $\mathcal{G}$ and $\mathcal{H}$ be marked graphs. The vertices $v_1,\dots, v_k$ of $\mathcal{G}$ are marked by the numbers $p_i^{n_i}$ and the vertices  $w_1,\dots, w_r$ of  $\mathcal{H}$ are marked by the numbers $q_j^{m_j}$, where $p_i, q_j$ are primes and $n_i,m_j\in \N$. We say that an embedding $\varphi:v_i\mapsto w_i$ of the graph $\mathcal{G}$ into the graph $\mathcal{H}$ is an \emph{embedding of marked graphs} if $q_i=p_i$ and $n_i\le m_i$ for all $i=1,\dots, k$. It is clear that an embedding of marked graphs induces a homomorphism from $\GG$ to $\HH$, where $\HH$ is the graph product defined by the marked graph $\mathcal{H}$.

If $\mathcal{G}'$ is a full subgraph of $\mathcal{G}$. Then the subgroup $\GG'$ of $\GG$ generated by the vertex groups of the graph $\mathcal{G}'$ is called a \emph{canonical parabolic subgroup} of $\GG$. One can show that $\GG'$ has the structure of a graph product of groups with the underlying graph $\mathcal{G}'$.

A word in $\GG$ is a finite string of elements from the alphabet consisting of the non-identity elements of $G_1,\dots, G_k$. It is clear that any word represents an element of $\GG$. A word $g_1 \dots g_n$ is called \emph{reduced} if and only if $g_i$ and $g_{i+1}$ are in different vertex groups for each $i = 1, \dots , n - 1$. A reduced word is called \emph{geodesic} if and only if it has minimum length among all reduced words for the element $w \in \GG$. We say that $w\in \GG$ is \emph{cyclically reduced} if and only if for any geodesic word $g_1 \cdots g_n$ for $w$, either $n = 1$ or $g_1 \ne g_{n}^{-1}$. Let $w\in \GG$ be an element of finite order. We denote the order of $w$ by $\ord(w)$.

Suppose that a non-trivial element $w \in \GG$ has a geodesic word $g_1 \dots g_n$. Then the \emph{alphabet of $w$}, $\az(w)$, is the set of vertex groups to which the $g_i$'s belong. It is not hard to show that if $w_1$, $w_2$ are minimal words that represent the same element in $\GG$, then $\az(w_1)=\az(w_2)$, see Proposition 2.4.2, \cite{Goda}.

Given the graph $\mathcal{G}$, we define \emph{the non-commutation graph} $\Gamma$ of $\GG$. The graph $\Gamma$ has the same vertex set as $\mathcal{G}$ but two vertices are joined by an edge in $\Gamma$ if and only if they are not joined by an edge in $\mathcal{G}$.

A cyclically reduced element $w$ of $\GG$ is called a \emph{block} if and only if the full subgraph on $\az(w)$ of the non-commutation graph is connected. A block is called \emph{regular} if and only if its length is greater than one. Equivalently, the block $w$ is regular if the full subgraph on $\az(w)$ has more than one vertex. Otherwise a block is called \emph{singular}. Denote by $\BA(w)$ the canonical parabolic subgroup of $\GG$ generated by all vertex groups that do not occur in $\az(w)$ and commute with $w$.

Every cyclically reduced element of $\GG$ can be written as a product of blocks in a unique way, see Section 2.10, \cite{Goda}.

We now give a description of centralisers of elements in $\GG$. The theorem below is a particular case of Theorem 5.6.4 in \cite{Goda} for graph products of finite abelian groups.

\begin{thm}[Centraliser Theorem]
Let $\GG$ be a graph product of finite abelian groups. Let $w$ be a cyclically reduced element of $\GG$ and suppose that $w = w^{(1)}\cdots w^{(t)}$ is its block decomposition. Index the blocks of $w$ so that $w^{(1)}, \dots , w^{(l)}$ for some $l\le t$ are the singular blocks of $w$ with each $w^{(j)}\in G_{i_j}$ and $w^{(k+1)},\dots , w^{(t)}$ are the regular blocks of $w$. Then the centraliser $C(w)$ of $w$ is
$$
C(w)=G_{i_1}\times \dots \times G_{i_l}\times \langle w^{(l+1)}\rangle \times \dots \times \langle w^{(t)}\rangle \times \BA(w).
$$
\end{thm}

The following proposition is a recapitulation of Theorem 3.26 in \cite{Green} for graph products of finite abelian groups. It is also an easy consequence of the above theorem.
\begin{prop} \label{prop:finord}
Let $\GG$ be the graph product of finite abelian groups. Let $g\in \GG$, $g=h^{-1}g'h$, where $g'$ is cyclically reduced. Then the order of $g$ is finite if and only if every block of $g'$ is singular. In other words, $g'$ belongs to a canonical parabolic subgroup of $\GG$, whose underlying graph is complete.
\end{prop}

\subsection{First-order logic}
\label{se:2-3}

In this section we recall some basic notions of first-order logic and model theory. We refer the reader to \cite{ChKe} for details.

The standard first-order language of group theory, which we denote by $\cL$, consists of a symbol for multiplication $\cdot$, a symbol for inversion $^{-1}$, and a symbol for the identity $1$.

Any group word in variables $X$ can be considered as as a term in the language $\cL$. Observe that every term in the language $\cL$ is  equivalent modulo the axioms of group theory to a group word in variables $X$. An {\em atomic formula}  in the language $\cL$ is a formula of the type $W(X) = 1$, where $W(X)$ is a group word in $X$.  A {\em Boolean combination} of atomic formulas in the language $\cL$ is a disjunction of conjunctions of atomic formulas and their negations. Thus every Boolean combination $\Phi$  of atomic  formulas in $\cL$ can be written in the form $\Phi =  \bigvee\limits_{i=1}^m\Psi_i$, where each $\Psi_i$ has one of following forms:
$$
\bigwedge\limits_{j = 1}^n(S_j(X) = 1),  \hbox{ or } \bigwedge\limits_{j =1}^n(T_j(X) \neq 1),  \hbox{ or }
\bigwedge\limits_{j = 1}^n (S_j(X) = 1) \ \wedge \  \bigwedge_{k = 1}^m (T_k (X) \neq 1).
$$

It follows from general results on disjunctive normal forms in propositional logic that every quantifier-free formula in the language $\cL$ is logically equivalent (modulo the axioms of group theory) to a Boolean combination of  atomic ones. Moreover, every formula $\Phi$ in $\cL$ with \emph{free variables} $Z=\{z_1,\ldots ,z_k\}$ is logically equivalent to a formula of the type
$$
Q_1x_1 Q_2 x_2 \ldots Q_l x_l \Psi(X,Z),
$$
where  $Q_i \in \{\forall, \exists \}$, and  $\Psi(X,Z)$ is a Boolean combination of atomic formulas in variables from $X \cup Z$. A first-order formula $\Phi$ is called a \emph{sentence}, if $\Phi$ does not contain free variables.

For a group $H$ the set of all sentences in $\cL$ which hold in $H$ is called the {\em elementary theory} of $H$ in the language $\cL$. Two groups $H$ and $K$ are called {\em elementarily equivalent } if and only if their elementary theories coincide.

A sentence $\Phi$ is called an \emph{$\exists\forall$-sentence} if and only if $\Phi$ has the form
$$
\Phi= \exists x_1 \dots \exists x_k\forall y_1 \dots\forall y_l \Psi(x_1,\dots, x_k, y_1, \dots, y_l),
$$
where $\Psi(x_1,\dots, x_k, y_1, \dots, y_l)$ is quantifier-free. The collection of all $\exists\forall$-sentences satisfied by a group $H$ is called the $\exists\forall$-theory of $H$. Two groups are called  \emph{$\exists\forall$-equivalent} if and only if their $\exists\forall$-theories coincide. If a sentence $\Phi$ holds in a group $H$, then we write $H \models \Phi$.

\section{Elementary equivalence of graph products}
Throughout this section let $\GG$ and $\HH$ be the graph product of finite abelian groups. Denote the corresponding marked graphs by $\mathcal{G}$ and $\mathcal{H}$. The marked graph $\mathcal{G}$ has vertices $v_1,\dots, v_k$, which are labelled by the numbers $p_i^{n_i}$, where $p_i$ is a prime and $n_i\in \N$. Similarly, the marked graph $\mathcal{H}$ has vertices $w_1,\dots, w_r$, which are labelled by the numbers $q_i^{m_i}$, where $q_i$ is a prime and $m_i\in \N$. Therefore, $\GG$ and $\HH$ are graph products of finite, cyclic, directly indecomposable groups $G_1,\dots, G_k$ and $H_1,\dots, H_r$, respectively.

By the graph product $\GG$ we write a sentence $\Phi_\GG=\Phi_\GG(x_1,\dots, x_k)$ that, as we shall see later, describes the group $\GG$ up to isomorphism. The formula $\Phi_\GG$ states that there exist $k$ elements $x_1, \dots, x_k$ (recall that the number of vertices of the marked graph $\mathcal{G}$ underlying the construction of the graph product is $k$) so that
\begin{enumerate}
    \item\label{it:1} the order of  $x_i$ equals  ${\ord(G_i)}$;
    \item\label{it:2} if there is an edge in $\mathcal{G}$ between the vertices $v_i$ and $v_j$, then $x_i$ and $x_j$ commute;
    \item\label{it:3} for all $g_1,\dots, g_{k-1}$, the element $x_i^{s}$ is not equal to ${\left(x_{i_1}^{t_{i_1}}\right)}^{g_{i_1}}\cdots {\left(x_{i_l}^{t_{i_l}}\right)}^{g_{i_l}}$, where $x_{i_{j_1}}\ne x_{i_{j_2}}$ if $j_1\ne j_2$, $i_j\ne i$ for all $1\le j,j_1,j_2\le l$, and $1\le t_{i_j}\le \ord(G_{i_j})-1$, $0\le s\le \ord(G_i)-1$.
\end{enumerate}
It is clear that the above can be written using the first-order language of groups. Furthermore, the formula $\Phi_\GG$ is an $\exists\forall$-formula. Note that if $g_i$ is a generator of the group $G_i$, then $\GG\models \Phi_\GG(g_1,\dots, g_k)$, in other words the tuple of elements $g_1,\dots, g_k$ satisfies conditions (\ref{it:1}),  (\ref{it:2}) and  (\ref{it:3}).

\begin{prop}\label{prop:main}
Let $\GG$ and $\HH$ be two graph products of finite abelian groups defined by marked graphs $\mathcal{G}$ and $\mathcal{H}$ respectively. If $\HH\models \Phi_\GG (a_1,\dots, a_k)$, then there exist singular block elements $h_1,\dots, h_k$ such that  $\HH\models \Phi_\GG (h_1,\dots, h_k)$.
\end{prop}
\begin{proof}
Suppose that $\HH\models \Phi_\GG (a_1,\dots, a_k)$. By condition (\ref{it:1}) from the construction of $\Phi_\GG$ we have that the order of $a_i$ equals $\ord(G_i)$. Therefore, by Proposition \ref{prop:finord}, we get that $a_i$ belongs to a conjugate $C_i^{f_i}$ of a canonical parabolic subgroup $C_i$ of $\HH$, whose underlying graph is complete, $1\le i\le k$, i.e. $a_i=c_i^{f_i}$, $c_i \in C_i$. We now show that $\HH\models \Phi_\GG (c_1,\dots, c_k)$. Indeed, since conditions (\ref{it:1}) and (\ref{it:3}) from construction of $\Phi_\GG$ are invariant under conjugation, the tuple $c_1,\dots, c_k$ satisfies conditions (\ref{it:1}) and (\ref{it:3}). Furthermore, by the Centraliser Theorem, if $a_i$ and $a_j$ commute, then $c_i$ and $c_j$ commute, therefore  the tuple $c_1,\dots, c_k$ satisfies condition (\ref{it:2}). Note that every block of $c_i$ is singular.

We now prove that $\HH\models \Phi_\GG(d_1,\dots, d_k)$, where the tuple of elements $d_1,\dots, d_k$ satisfies that the order of any block of $d_i$ equals the order of $d_i$. Let
$$
c_i=b_{i,1}\cdots b_{i,l_i}
$$
be the block decomposition of $c_i$. Since $\ord(c_i)=p_i^{n_i}$ and $p_i$ is a prime, we have that $\ord(b_{i,j})=p_i^{n'_i}$, where $0<n'_i\le n_i$. Furthermore, there exists a block $b_{i,j_0}$ of $c_i$ such that $\ord(b_{i,j_0})=p_i^{n_i}$. Set $d_i$ to be the product of the blocks of $c_i$ whose order equals $p_i^{n_i}$. We now show that $\HH\models \Phi_\GG(c_1,\dots, c_{i-1},d_i,c_{i+1},\dots, c_k)$.  It is clear that the tuple $c_1,\dots, c_{i-1},d_i,c_{i+1},\dots, c_k$ satisfies conditions (\ref{it:1}) and (\ref{it:2}) from the definition of $\Phi_\GG$. Suppose that condition (\ref{it:3}) fails. Without loss of generality, we may assume that
$$
d_i^{s}={\left(c_{i_1}^{t_{i_1}}\right)}^{g_{i_1}}\cdots {\left(c_{i_l}^{t_{i_l}}\right)}^{g_{i_l}}.
$$
Write $s=s_1p_i^{s_2}$, where $s_1$ and $p_i$ are co-prime (note that $s_2<n_i$). Let
$$
M_i=\max\left\{\max\limits_{1\le j\le l_i}\left\{\ord(b_{i,j})\mid \ord(b_{i,j})<p_i^{n_i}\right\}, p_i^{s_2}\right\}.
$$
We have that $s_1M_i\ne 0\mod(p_i^{n_i})$. Therefore,
$$
1\ne d_i^{T}=d_i^{s_1M_i}={\left({\left(c_{i_1}^{t_{i_1}}\right)}^{g_{i_1}}\cdots {\left(c_{i_l}^{t_{i_l}}\right)}^{g_{i_l}}\right)}^{M_i'}= {\left(c_{i_1}^{\bar t_{i_1}}\right)}^{g_{i_1}}\cdots {\left(c_{i_l}^{\bar t_{i_l}}\right)}^{g_{i_l}},
$$
where $T=s_1M_i\mod (\ord(G_i))$, $M_i'=\frac{M_i}{p_i^{s_2}}$ and $\bar t_{i_j}=t_{i_j}M_i'\mod(\ord(G_{i_j}))$. Since $d_i^{T}=c_i^T$, we obtain a contradiction with the assumption that $\HH\models \Phi_\GG(c_1,\dots, c_k)$. Therefore, $\HH\models \Phi_\GG(c_1,\dots, c_{i-1},d_i,c_{i+1},\dots, c_k)$. Recursively applying the above argument we get that $\HH\models \Phi_\GG(d_1,\dots, d_k)$.

Finally, we show that there exists a tuple of singular block elements $h_1,\dots, h_k$ such that $\HH\models \Phi_\GG(h_1,\dots, h_k)$. Fix generators $u_1,\dots,u_r$ of the finite, cyclic, directly indecomposable groups $H_1,\dots, H_r$.  Consider the block decomposition of $d_i=b_{i,1}\dots b_{i,l_i}$. Then, since $b_{i,j}$ is a singular block, we have  $b_{i,j}=u_j^{l_{i,j}}$, where $1\le l_{i,j}\le \ord(H_j)-1$.

Consider a $k\times r$-matrix $M$, where $k$ and $r$ are the numbers of vertices of the marked graphs $\mathcal{G}$ and $\mathcal{H}$, respectively. The element $e_{i,j}$ of $M$ equals $l_{i,j}$ if and only if there exists a block $b_{i,j}=u_j^{l_{i,j}}$ in the block decomposition of $d_i$. Otherwise set $e_{i,j}=0$. Since the tuple $d_1,\dots, d_k$ satisfies condition  (\ref{it:3}) from the construction of the formula $\Phi_\GG$, it follows that the rows of the matrix $M$ are linearly independent. Therefore, the rank of $M$ equals $k$.

We use induction to show that one can choose a tuple of non-zero elements $e_{1,j_1}, \dots, e_{k,j_k}$ such that $j_{i_1}\ne j_{i_2}$. The case $k=1$ is trivial. Choose $k$ columns of $M$ so that the $k\times k$ matrix $N_k$ consisting of these $k$ columns has rank $k$. The determinant  of $N_k$ equals $e_{1,1} L_1 -e_{1,2}L_2+\dots +(-1)^{k+1} e_{1,k}L_k$, where $L_i$ is the corresponding minor. Since the determinant of $N_k$ is non-zero,  there exists $j$ such that $e_{1,j}\ne 0$ and $L_j\ne 0$. Since $L_j$ is the determinant of a $(k-1)\times (k-1)$ matrix, by induction, we can choose $k-1$ non-zero elements $e_{2,j_2},\dots, e_{k,j_k}$ such that $j_{i_1}\ne j_{i_2}$, $j_{i_1}, j_{i_2}\ne j$ and the statement follows.

Therefore, if we set $h_1=b_{1,{j_1}}, \dots, h_k=b_{k,j_k}$, then the tuple of elements $h_1,\dots, h_k$ satisfies conditions (\ref{it:1})-(\ref{it:3}) from the definition of $\Phi_\GG$.
\end{proof}

\begin{cor} \label{cor:main}
Let $\GG$ and $\HH$ be two graph products of finite abelian groups defined by marked graphs $\mathcal{G}$ and $\mathcal{H}$ respectively. If $\HH\models \Phi_\GG (a_1,\dots, a_k)$, then the marked graph $\mathcal{G}$ embeds into the marked graph $\mathcal{H}$.
\end{cor}

\begin{thm}
Let $\GG$ and $\HH$ be two graph products of finite abelian groups with the underlying marked graphs $\mathcal{G}$ and $\mathcal{H}$ respectively. The following are equivalent:
\begin{itemize}
\item the marked graphs $\mathcal{G}$ and $\mathcal{H}$ are isomorphic;
\item the groups $\GG$ and $\HH$ are isomorphic;
\item the groups $\GG$ and $\HH$ are elementarily equivalent;
\item the groups $\GG$ and $\HH$ are $\exists\forall$-equivalent.
\end{itemize}
\end{thm}
\begin{proof}
If the marked graphs $\mathcal{G}$ and $\mathcal{H}$ are isomorphic, then it is clear that the groups $\GG$ and $\HH$ are isomorphic. If two groups are isomorphic, then they are elementarily equivalent. If two groups are elementarily equivalent, then, in particular, they are $\exists\forall$-equivalent.

We now show that if $\GG$ and $\HH$ are $\exists\forall$-equivalent, then the marked graphs $\mathcal{G}$ and $\mathcal{H}$ are isomorphic. Since $\Phi_\GG$ and $\Phi_\HH$ are $\exists\forall$-formulas, if the groups $\GG$ and $\HH$ are $\exists\forall$-equivalent, then $\GG\models \Phi_\HH$ and $\HH\models \Phi_\GG$. By Corollary \ref{cor:main}, it follows that the marked graph $\mathcal{G}$ embeds into $\mathcal{H}$ and vice-versa, hence $\mathcal{G}$ and $\mathcal{H}$ are isomorphic.
\end{proof}

In particular, the above theorem shows that two graph products of finite abelian groups are isomorphic if and only if the corresponding marked graphs are. A generalisation of this result to graph products of finite directly indecomposable groups was proved in \cite{Radcliffe}.

\begin{rem}
The authors believe that the main result of this paper holds for graph products of arbitrary finite groups. The argument, however, would be more technical. In particular, condition (\ref{it:1}) in the construction of the formula $\Phi_\GG$ should describe the multiplication table of the corresponding finite vertex group.
\end{rem}

\bibliographystyle{amsplain}

\end{document}